\documentclass[oneside,english]{amsart}
\usepackage[T1]{fontenc}
\usepackage[latin1]{inputenc}
\usepackage{a4wide}
\usepackage{amssymb}

\makeatletter

 \theoremstyle{plain}
 \theoremstyle{plain}    
 \newtheorem*{thm*}{Theorem} 

\usepackage{babel}
\makeatother
\begin{document}

\title{On the origin of the BV operator on odd symplectic supermanifolds}

\author{Pavol \v Severa}

\address{Dept.~of theoretical physics \\
Mlynská dolina F2, 84248 Bratislava\\
Slovakia}

\email{severa@sophia.dtp.fmph.uniba.sk}

\maketitle
The invariant meaning of the Batalin-Vilkoviski operator was discovered
by H.~Khudaverdian \cite{kh}. He noticed that if $x^{i}$, $p_{i}$
are local Darboux coordinates on a supermanifold $M$ with an odd
symplectic structure, the operator\[
\Delta=\frac{\partial^{2}}{\partial x^{i}\partial p_{i}},\]
acting on \emph{semidensities} on $M$, is independent of the choice
of coorditates. The aim of this note is to find a natural definition
of this somewhat miraculous operator.

\subsection*{Two differentials on odd symplectic supermanifolds}

In what follows, let $M$ be a supermanifold with an odd symplectic
form $\omega$. On $\Omega(M)$ we have two anticommuting differentials:
one is de Rham's $d$ and the other is the wedge product with $\omega$.
We shall find that the Khudaverdian BV operator is (rougly speaking)
the 3rd differential of the spectral sequence of this bicomplex. 

\begin{thm*}
Let $(M,\omega)$ be an odd symplectic supermanifold. In the spectral
sequence of the bicomplex $\left(\Omega(M),\,\omega\wedge,\, d\right)$
we have:
\begin{enumerate}
\item the cohomology of the complex $\left(\Omega(M),\,\omega\wedge\right)$
is naturally isomorphic to the semidensities on $M$
\item the next differential in the spectral sequence, de Rham's $d$, vanishes
on the cohomology of $\left(\Omega(M),\,\omega\wedge\right)$
\item the next differential ($d\circ\left(\omega\wedge\right)^{-1}\circ d$)
coincides with the BV operator
\item all higher differentials are zero.
\end{enumerate}
\end{thm*}
The proof of this theorem is completely straightforward; we shall
do it leisurely in the rest of this note.

\subsection*{Cohomology of $\omega\wedge$}

It is fairly simple to describe the cohomology of the complex $\left(\Omega(M),\,\omega\wedge\right)$
in local Darboux coordinates $x^{i}$, $p_{i}$ ($i=1,\dots,n$, $\omega=dp_{i}\wedge dx^{i}$),
where $x^{i}$ are the even coordinates and $p_{i}$ the odd coordinates.
Let $U\subset M$ be the open subset covered by the coordinates. Then
any cohomology class of $\left(\Omega(U),\,\omega\wedge\right)$ has
unique representative of the form \begin{equation}
f(x,p)\; dx^{1}\wedge dx^{2}\wedge\dots\wedge dx^{n}.\label{eq:rep}\end{equation}
In other words, using the coordinates, we can locally identify $H\left(\Omega(M),\,\omega\wedge\right)$
with functions on $M$ (this identification \emph{does} depend on
the choice of coordinates). 

To prove this claim we split $\left(\Omega(U),\,\omega\wedge\right)$
into subcomplexes: we assign an auxiliary degree $1$ to each $dx$
and $-1$ to each $dp$, and denote this degree $\mathrm{auxdeg}$
($x$'s and $p$'s will have $\mathrm{auxdeg}=0$); the subspaces
of $\Omega(U)$ of fixed $\mathrm{auxdeg}$ are clearly subcomplexes.
We shall see that each of them has zero cohomology, except for the
one with degree $n$, where the differential vanishes. We shall prove
it using an explicit homotopy. Let us consider the operator $L:\Omega(U)\rightarrow\Omega(U)$
given by\[
L:\alpha\mapsto\partial_{x^{i}}\lrcorner\,\partial_{p_{i}}\lrcorner\,\alpha.\]
A direct computation shows that\[
L\circ(\omega\wedge)+(\omega\wedge)\circ L:\alpha\mapsto(n-\mathrm{auxdeg}\,\alpha)\,\alpha.\]
This concludes the proof.

Now we also see that $d$ is 0 on $H\left(\Omega(M),\,\omega\wedge\right)$,
since \[
d(f(x,p)\; dx^{1}\wedge dx^{2}\wedge\dots\wedge dx^{n})=\frac{\partial f}{\partial p_{k}}\; dp_{k}\wedge dx^{1}\wedge dx^{2}\wedge\dots\wedge dx^{n},\]
which is $\omega\wedge$-exact (having $\mathrm{auxdeg}=n-1$).

\subsection*{The third differential}

Let us now compute the 3rd differential $d\circ\left(\omega\wedge\right)^{-1}\circ d$
in local Darboux coordinates. We have\[
d(f(x,p)\; dx^{1}\wedge dx^{2}\wedge\dots\wedge dx^{n})=\omega\wedge\alpha,\]
where\[
\alpha=L(d(f(x,p)\; dx^{1}\wedge dx^{2}\wedge\dots\wedge dx^{n}))=\frac{\partial f}{\partial p_{k}}\;\partial_{x^{k}}\lrcorner\, dx^{1}\wedge dx^{2}\wedge\dots\wedge dx^{n}\]
and $d\alpha$ is (up to a $\omega\wedge$-exact term)\[
\frac{\partial^{2}f}{\partial x^{k}\partial p_{k}}\; dx^{1}\wedge dx^{2}\wedge\dots\wedge dx^{n}.\]
The third differential in the spectral sequence is thus equal to the
Batalin-Vilkoviski operator\[
\Delta=\frac{\partial^{2}}{\partial x^{k}\partial p_{k}}.\]

Now if $M$ is contractible and admits global Darboux coordinates,
the cohomology of $\Delta$ is isomorphic to $\mathbb{R}$ (since
$\Delta$ can be identified with de Rham's $d$ on a contractible
subset of $\mathbb{R}^{n}$), and any cohomology class has a representative
in $\Omega(M)$ which is a constant multiple of $dx^{1}\wedge dx^{2}\wedge\dots\wedge dx^{n}$.
This representative is $d$-closed and thus is annuled by all higher
differentials in the spectral sequence. Since any $M$ can be covered
by such patches, these higher differentials vanish for any $M$. This
concludes the proof of the theorem, except for the part (1).

\subsection*{Semidensities}

Now we'll prove the part (1) of the theorem. We locally identified
the cohomology of $H\left(\Omega(M),\,\omega\wedge\right)$ with functions
on $M$ by choosing the representant (\ref{eq:rep}); it is fairly
easy to see that when we pass to another system of local Darboux coordinates,
the function $f$ gets multiplied by the square root of the corresponding
Berezinian. We shall however give a different proof, using Manin's
cohomological definition of Berezinian \cite{ma}. The claim we are
proving here, together with the proof, is taken from \cite{se}.

Let us recall Manin's definition. Let $V$ be a vector supespace.
Let us choose a vector superspace $W$ with an odd symplectic form
$\omega\in\bigwedge^{2}W^{*}$, such that $V$ is its Lagrangian subspace
(for example $W=V\oplus\Pi V^{*}$). Then $\mathrm{Ber}(V^{*})$ (the
1-dimensional vector space of constant densities on $V$) is defined
as the cohomology $H\left(\bigwedge W^{*},\,\omega\wedge\right)$
(this definition is easily seen to be independent of the choice of
$W$).

If now $V'$ is a Lagrangian complement of $V$ in $W$, then again
$\mathrm{Ber}(V'^{*})=H\left(\bigwedge W^{*},\,\omega\wedge\right)$;
on the other hand, \[
\mathrm{Ber}(W^{*})=\mathrm{Ber}(V^{*})\otimes\mathrm{Ber}(V'^{*})=H\left({\textstyle \bigwedge}W^{*},\,\omega\wedge\right)^{\otimes2},\]
and thus $H\left(\bigwedge W^{*},\,\omega\wedge\right)=\mathrm{Ber}(W^{*})^{1/2}$
(we should write everywhere {}``naturally isomorphic'' instead of
{}``equal'', but hopefully it's not a big crime). This identity
is valid for any vector superspace with an odd symplectic form. We
apply it to the bundle of symplectic vector spaces $TM$, which concludes
the proof.

\subsection*{Final remarks}

1. We should say a few remarks about the spectral sequence of the
bicomplex $\left(\Omega(M),\,\omega\wedge,\, d\right)$, since $\Omega(M)$
is \emph{not} bigraded. The spectral sequence is constructed in this
way: we take $\Omega(M)[\hbar]$ (differential forms on $M$ depending
polynomially on an indeterminate $\hbar$; we could just as well take
$\Omega(M)[\![\hbar]\!]$), with the differential $\hbar d+\omega\wedge$.
Then multiplication by $\hbar$ is an endomorphism of the complex
$(\Omega(M)[\hbar],\hbar d+\omega\wedge)$, and our spectral sequence
is the Bochstein spectral sequence of this endomorphism. That is,
we start with the short exact sequence of complexes\[
0\rightarrow(\Omega(M)[\hbar],\hbar d+\omega\wedge)\stackrel{\hbar\cdot}{\rightarrow}(\Omega(M)[\hbar],\hbar d+\omega\wedge)\rightarrow(\Omega(M),\omega\wedge)\rightarrow0,\]
 out of which we get the exact couple\[
\begin{array}{ccccc}
H(\Omega(M)[\hbar],\hbar d+\omega\wedge) &  & \stackrel{(\hbar\cdot)_{*}}{\longrightarrow} &  & H(\Omega(M)[\hbar],\hbar d+\omega\wedge)\\
 & \nwarrow &  & \swarrow\\
 &  & H(\Omega(M),\omega\wedge)\end{array}\]
which generates the spectral sequence. If we denote $E_{\infty}$
its ultimate term, we have \[
H\left(\Omega(M)[\hbar,\hbar^{-1}],\hbar d+\omega\wedge\right)\cong E_{\infty}[\hbar,\hbar^{-1}]\]
 (and also $H\left(\Omega(M)[\![\hbar]\!][\hbar^{-1}],\hbar d+\omega\wedge\right)\cong E_{\infty}[\![\hbar]\!][\hbar^{-1}]$).

2. The odd symplectic form $\omega$ on $M$ gives us an isomorphism
between $\Omega(M)$ and $\Gamma(STM)$, i.e.~the space of polynomial
functions on $T^{*}M$. This isomorphism transfers $\omega\wedge$
to multiplication by the odd Poisson structure $\pi$ corresponding
to $\omega$ (recall that since $\pi$ is an odd Poisson structure,
it is a function on $T^{*}M$), and $d$ to $\{\pi,\cdot\}$ (where
$\{,\}$ is the Poisson bracket on $T^{*}M$); the differential $\hbar d+\omega\wedge$
becomes $\pi+\hbar\{\pi,\cdot\}$. This suggests some generalizations,
e.g.~we can take an odd Poisson structure which is not symplectic,
or more generally, instead of $T^{*}M$ we can take an aritrary supermanifold
with an even symplectic form, and consider on it an odd function $\pi$
such that $\{\pi,\pi\}=0$. The spectral sequences would still be
defined, but it is not clear to me if they are good for anything.

3. The result in this note is extremely simple; it is written with
the hope that it might be helpfull in situations where BV-like operators
are much less trivial.

\end{document}